\documentclass[11pt,reqno]{article}
\usepackage{amssymb, amsmath, amsthm, amsfonts, amscd}
\usepackage[english]{babel}

\numberwithin{equation}{section}
\newtheorem{remark}{Remark}

\def\N{\mathbb{N}}

\def\R{\mathbb{R}}

\def\de{\delta}

\newcommand{\p}{\partial}
\newcommand{\pd}[2]{\frac {\p #1}{\p #2}}
\newcommand{\ds}{\displaystyle}

\newcommand{\Lcal}{\mathcal{L}}
\newcommand{\Ical}{\mathcal{I}}
\newcommand{\Fourier}{\mathcal{F}}
\newcommand{\Pcal}{\mathcal{P}}
\newcommand{\Qcal}{\mathcal{Q}}
\newcommand{\tLa}{\widetilde\Lcal_a}

\newcommand{\ga}{\gamma}

\newcommand{\Id}{\mathcal{I}}
\newcommand{\imaginary}{{\tt i}}

\numberwithin{equation}{section}

\begin{document}
\title{Photoacoustic imaging in attenuating acoustic media based on strongly causal models}
\author{Konstantinos Kalimeris \thanks{Radon Institute of Computational and Applied Mathematics, Altenberger Str.~69, 4040 Linz, Austria}
        \and Otmar Scherzer\thanks{Computational Science Center, University of Vienna, Nordbergstr.~15, 1090 Vienna, Austria and 
             Radon Institute of Computational and Applied Mathematics, Altenberger Str.~69, 4040 Linz, Austria}
        }
\date{\today}
\maketitle

\begin{abstract}
In this paper we derive time reversal imaging functionals for two \emph{strongly causal} acoustic attenuation models, which have been proposed recently. 
The time reversal techniques are based on recently proposed ideas of Ammari et al for the thermo-viscous wave equation. Here and there an asymptotic analysis 
provides reconstruction functionals from first order corrections for the attenuating effect. In addition, we present a novel approach for higher order 
corrections.
\end{abstract}

\section{Introduction}

\emph{Photoacoustic Imaging} is a promising imaging technique for visualizing biological material parameters. In
experiments, the medium is exposed to a short pulse of an electromagnetic wave. The medium absorbs a fraction of 
the induced energy, heats up, and reacts with thermoelastic expansion. This in turn produces acoustic waves, which 
can be recorded and which are used to determine the electromagnetic absorption coefficient 
\cite{XuWan06,ElbSchSchu11_report,KirSch11_report,ZanSchHal09b}. 
These coupling properties explain why photoacoustic imaging is referred to as \emph{hybrid} or \emph{Coupled Physics Imaging}. 
For some recent progress in hybrid imaging we refer to the surveys \cite{Kuc11,ArrSch12}.

In this paper we investigate the method of time reversal in attenuating media as it was introduced in 
\cite{AmmBreGarWah11,Wah11,AmmBreGarWah12_report} for the thermo-viscous wave equation. 
In these references, the goal is to construct a parameter dependent family of approximate reconstruction functionals, 
which allow for approximation of the initial datum of the thermo-viscous wave equation. 
For time-reversal technique for imaging in general, see for example \cite{Amm08,Fin97,FouGarPapSol07}. 

In this work, we are using the asymptotic techniques from \cite{AmmBreGarWah11,Wah11} to develop time reversal algorithms for the 
Nachman-Smith-Waag (NSW) \cite{NacSmiWaa90} and the Kowar-Scherzer-Bonnefond (KSB) models \cite{KowSch10,KowSchBon11}. 
These two models satisfy a \emph{strong causality property}. That is, the solutions of these equations are zero before the initialization and satisfy 
\emph{a finite front wave propagation speed}. We emphasize that the thermo-viscous wave equation considered in \cite{AmmBreGarWah11,Wah11} 
is not strongly causal \cite{IgnOst10}.
While there is a large literature on inversion formulas and time reversal algorithms for the standard wave equation much less is 
known in the case of attenuating media \cite{BurGruHalNusPal07,KowSch10,AmmBreGarWah11,NacSmiWaa90}. Partially, this is due to the 
fact that, so far, the reference attenuation model (that is, the governing wave equation) has not been established in this field. 

\section{Review on time-reversal for the thermo-viscous equation}
 
Let $\Omega$ be a bounded domain in $\R^3$. Standard photoacoustic imaging consists in determining the \emph{absorption density} $f$, which is assumed to have compact support $K$ in $\Omega$, 
in the acoustic wave equation
\begin{equation}\label{wave-eq}\begin{split}
\frac{\p^2 p}{\p t^2}(x,t) - \Delta p(x,t) = \pd{\de_0}{t}(t) f(x), \qquad (x,t)\in\R^3\times\R,\\
p(x,t)=0 \ \text{ and } \ \pd{p}{t}(x,t)=0,\qquad x\in\R^3, \ t < 0,
\end{split}\end{equation}
from measurement data $g(y,t):=p(y,t)$ for some $y\in\partial \Omega$ and $t\in[0,T]$, where $T$ is supposed to be sufficiently large.

In the following, let $a:=a(x)$ be a positive function, describing \emph{attenuation}. 
The imaging method considered in \cite{AmmBreGarWah11,Wah11} consist in reconstructing $f$ from data
\begin{equation*}
g_a(y,t):=p^a(y,t) \text{ for all } y\in\partial \Omega \text{ and } t\in[0,+\infty).
\end{equation*}
where $p^a$ solves the \emph{thermo-viscous wave} equation,
\begin{equation} 
\label{wave-eq-att-tv} 
\begin{split}
\frac{\p^2 p^a}{\p t^2}(x,t) -\left(\Id + a \frac{\p}{\p t} \right) \Delta p^a(x,t) =   \pd{\delta_0}{t}(t) f(x),\qquad (x,t)\in\R^3\times\R,\\
p^a(x,t)=0 \ \text{ and } \ \pd{p^a}{t}(x,t)=0,\qquad x\in\R^3, \ t \ll 0\;.
\end{split}
\end{equation} 
To derive the imaging technique we use $\tilde\Gamma_\omega^a(x,\cdot)$, the fundamental solution at $x$
of the Helmholtz equation,
\begin{equation}
\label{eq:helmholtz}
\omega^2 \tilde u (x,y)+(1+\imaginary a \omega)\Delta_y \tilde u(x,y)= -  \delta_x(y), \qquad y\in \R^3\;.
\end{equation}
In \cite{AmmBreGarWah11,Wah11} the following results have been shown:
\begin{itemize}
\item For fixed cut off parameter $\rho > 0$ the function
\begin{equation} \label{rel-1}
v_{s,\rho}^a(x,t):=-\frac{1}{2\pi}\int_{-\rho}^{\rho}\int_{\partial \Omega} \imaginary \omega \tilde{\Gamma}_\omega^a(x,y)
g_a(y,T-s) d\sigma(y)  e^{-\imaginary \omega(t-s)} d\omega\,,
\end{equation}
satisfies the thermo-viscous wave equation
\begin{equation} \label{reg-treq-tv}
\frac{\p^2 v}{\p t^2}(x,t) -\left(\Id - a \frac{\p}{\p t} \right) \Delta v(x,t) =   S_\rho\left[\pd{\delta_s}{t}\right]
g_a(x,T-s)\delta_{\partial \Omega}\,, 
\end{equation}
where $S_\rho$ is as in \eqref{S_r}. 
In this formula the necessity of the regularization becomes evident because in the unregularized 
form the right hand side consists of a product of two distributions, which a-priori is not well-defined. 
\item Moreover, for small $a$ it follows from the results in \cite{Wah11} that
\begin{equation} \label{functional-reg}
\Ical_\rho^a(x)=\int_0^T v_{s,\rho}^a(x,T)ds \longrightarrow f(x), \qquad \text{as } \rho\to+\infty\,.
\end{equation}
Moreover, \cite[Remark 2.3.6]{Wah11} also gives a reference how to obtain appropriate values for the cut off parameter
$\rho$. The regularization is required so that the function $v_{s,\rho}^a(x,t)$ is well-defined in \eqref{rel-1}. 
\end{itemize}


\section{The KSB-model} \label{section-KSB}
In this section, we are deriving an imaging functional for the KSB model \cite[Eqs.(28),(88)]{KowSch10} following the time reversal 
approach of \cite{AmmBreGarWah11} for the thermo-viscous wave equation, as outlined above. 

Let $\alpha_0 > 0$ and $\ga\in(1,2]$ be fixed. Then the KSB model assumes that the attenuated pressure $p^a$ satisfies the equation
\begin{equation}
\label{wave-KSB}\begin{split} 
\left( \alpha_0 \Id + L^{1/2}\right)^2\frac{\p^2 p^a}{\p t^2}(x,t) -L \ \Delta p^a(x,t)=L\
\pd{\delta_0(t)}{t} f(x),\qquad x\in\R^3, \ t\in \R,\\
p^a(x,t)=0 \text{ and } \ \pd{p^a}{t}(x,t)=0,\qquad x\in\R^3, \ t<0,
\end{split} \end{equation}
where $L^{1/2}$ is the convolution operator (in time) with kernel $\ds \frac{1}{\sqrt{2\pi}} \Fourier^{-1}\left[
\left(1+(-\imaginary \tau_0 \omega)^{\gamma-1}\right)^{1/2}\right]$ - we emphasize that $L=L^{1/2} \circ L^{1/2}$. 
Here $D_t^{\gamma-1}$ denotes a fractional time derivative operator of order $\gamma-1$ \cite{KilSriTru06,Pod99}. 

The Fourier transforms, $\hat{p}^a:=\Fourier[p^a](\omega)$ and $\hat{p}:=\Fourier[p](\omega)$, of 
\eqref{wave-KSB} and \eqref{wave-eq} satisfy the Helmholtz equations:
\begin{equation*}
\kappa^2(\omega) \hat{p}^a(x)+\Delta \hat{p}^a(x)= \imaginary \omega f(x), \text{ and } \omega^2
\hat{p}(x)+\Delta \hat{p}(x)= \imaginary \omega f(x)\,,
\end{equation*}
respectively. Here
\begin{equation}
\label{eq:kappa}
\kappa(\omega)=\omega \left(1+\frac{\alpha_0}{\left(1+(-\imaginary \tau_0 \omega)^{\gamma-1}\right)^{1/2}} \right)\;.
\end{equation}
Using the particular form of the Helmholtz equations it follows that 
\begin{equation*}
\hat{p}^a=\frac{\omega}{\kappa(\omega)}\Fourier[p](\kappa(\omega))\,,
\end{equation*}
which after applying the inverse Fourier transform yields
\begin{equation*}
p^a(x,t)=\Lcal_a[p(x,\cdot)](t)
\end{equation*}
with
\begin{equation} 
\label{operator1}
\Lcal_a[\phi](t) :=\Fourier^{-1} \left[ \frac{\omega}{\kappa(\omega)} \Fourier[\phi](\kappa(\omega))\right](t)
=\frac{1}{2\pi} \int_\R \frac{\omega}{\kappa(\omega)}  e^{-\imaginary \omega t} \int_\R e^{\imaginary
\kappa(\omega)s}\phi(s)dsd\omega\;.
\end{equation}
The goal is to find an imaging operator of the form, that is functions $\tilde{\kappa}$ and $\lambda$, 
\begin{equation} \label{tilde-La-operator}
\tLa[\phi](t) 
:= \frac{1}{2\pi}
\int_\R \frac{\omega \lambda(\omega)}{\tilde{\kappa}(\omega)} e^{-\imaginary \omega t} \int_\R \phi(s)e^{\imaginary
\tilde{\kappa}(\omega)s} ds d\omega
=
\Fourier^{-1} \left[ \frac{\omega \lambda(\omega)  }{\tilde{\kappa}(\omega)}
\Fourier[\phi](\tilde{\kappa}(\omega))\right](t)\,,
\end{equation}
for which the following expansion with respect to $\alpha_0$ (as in \eqref{eq:kappa}) holds:
\begin{equation}
\label{id1}
{\tLa}^* \Lcal_a[\phi](t)= \phi(t) + o(\alpha_0), \quad as \ \alpha_0 \to 0\;.
\end{equation}
Thereby, 
\begin{equation}
\label{operator2}
{\tLa}^*[\phi](t) :=
\frac{1}{2\pi} \int_\R \frac{\omega}{\tilde{\kappa}(\omega)}\lambda(\omega)  e^{\imaginary \tilde{\kappa}(\omega) t}
\int_\R e^{-\imaginary \omega s}\phi(s)dsd\omega
\end{equation}
is the adjoint of $\tLa$.

The derivation of imaging operators as in \eqref{tilde-La-operator} follows general principles, which are used later on for the 
other attenuation models as well. 
The principle consists in construction auxiliary functions $\lambda_1$, $\nu_1$, and $\nu_2$ which are determined from 
$\kappa$ by general construction without taking into account the special structure of the function. The explicit construction 
comes at a later stage.
We introduce the auxiliary function $\lambda_1: \R \to \R$ which is related to $\kappa$ in the following way:
\begin{equation}
\label{order}
\kappa(\omega)=\omega(1-\alpha_0 \lambda_1(\omega))+{\mathcal O}(\alpha_0^2).
\end{equation}
Then, it follows by expansion with respect to $\alpha_0$ that 
\begin{equation*}
\begin{aligned}
\frac{\omega}{\kappa(\omega)} &= 1 + \alpha_0 \lambda_1(\omega) + {\mathcal O}(\alpha_0^2) \text{ for fixed } \omega \in
\R\,,\\
e^{\imaginary \kappa(\omega) s} &= e^{\imaginary \omega s} (1 - \imaginary \alpha_0 \omega \lambda_1(\omega) s) +
{\mathcal O}(\alpha_0^2)  \text{ for fixed } s \in \R\,,
\end{aligned}
\end{equation*}
and consequently for fixed $s,\omega \in \R$
\begin{equation}
\label{prod1}
\frac{\omega}{\kappa(\omega)} e^{\imaginary \kappa(\omega)s} = e^{\imaginary \omega s} \left(
1 + \alpha_0 (\lambda_1(\omega) - \imaginary \omega \lambda_1(\omega) s)\right) + {\mathcal O}(\alpha_0^2)\;.
\end{equation}
This shows that
\begin{equation*}
\begin{aligned}
\Lcal_a[\phi(s)](t) &= \phi(t) + \alpha_0  \Fourier^{-1}\left[ \lambda_1(\omega) \Fourier[\phi(s)](\omega)-\imaginary
					 \omega\lambda_1(\omega) \Fourier[s\phi(s)](\omega)\right](t) +o(\alpha_0)\\
                    &=\Lcal_0[\phi(s)](t) +
					 \alpha_0  \Fourier^{-1}\left[ \lambda_1(\omega) \Fourier[\phi(s)](\omega)-\imaginary
					 \omega\lambda_1(\omega) \Fourier[s\phi(s)](\omega)\right](t)
					 + o(\alpha_0)\;.
\end{aligned}
\end{equation*}
Note that for $\alpha=0$, $\Lcal_0 = \Id$.
\footnote{
In \cite{AmmBreGarWah11,Wah11}, $\big(\kappa(\omega)-\omega\big)$ is an imaginary function, therefore the Stationary Phase Method is used for the asymptotic analysis of the previous integral operators. 
Here the relevant function, $\big( -\alpha_0 \lambda_1(\omega) \big)$ in its general form, is complex (with non-vanishing real part) and the usage of the method of Steepest Descent would provide the relevant 
asymptotic expansion. The Stationary Phase Method can be seen as a special case of the method of Steepest Descent. In this work, the relevant real part is zero for the NSW and the thermo-viscous model and consequently the Stationary Phase Method can be applied. However for the KSB model this part is non-zero and the use of the method of Steepest Descent is needed.}
In order to get an explicit form of $\tLa$ we introduce auxiliary functions $\nu_i: \R \to \R$, $i=1,2$, and 
\begin{equation}
\label{tkappa}
\begin{aligned}
\tilde{\kappa}(\omega) &:=\omega(1-\alpha_0 \nu_1(\omega))+{\mathcal O}(\alpha_0^2)\,,\\
\lambda(\omega)&:=1+\alpha_0 \nu_2(\omega)+{\mathcal O}(\alpha_0^2).
\end{aligned}
\end{equation}
Again, by expansion with respect to $\alpha_0$ it follows that
\begin{equation*}
\begin{aligned}
e^{\imaginary \tilde{\kappa}(\omega) s}&=e^{\imaginary \omega s} (1+\alpha_0 (-\imaginary \omega\nu_1(\omega))
s)+{\mathcal O}(\alpha_0^2)\,,\\
\frac{\omega\lambda(\omega)}{\tilde{\kappa}(\omega)}&=1+\alpha_0 (\nu_1(\omega)+\nu_2(\omega))+{\mathcal
O}(\alpha_0^2)\,,
\end{aligned}
\end{equation*}
and consequently,
\begin{equation}
\label{nu}
\frac{\omega\lambda(\omega)}{\tilde{\kappa}(\omega)}e^{\imaginary \tilde{\kappa}(\omega) s}
=e^{\imaginary \omega s} (1+\alpha_0(\nu_1(\omega)+\nu_2(\omega)-\imaginary \omega \nu_1(\omega)s))+{\mathcal
O}(\alpha_0^2)\;.
\end{equation}
Therefore, we have
\begin{equation*}
\begin{aligned}
{\tLa}^*[\phi(s)](t) =& \phi(t) \\
&+ \alpha_0 \left\lbrace \Fourier^{-1}\left[ (\nu_1(-\omega)+\nu_2(-\omega)) \Fourier[\phi(s)](\omega)\right] (t)
\right.\\
&\qquad + t
\left. \Fourier^{-1}\left[\imaginary \omega \nu_1(-\omega) \Fourier[\phi(s)](\omega)\right](t)  \right\rbrace \\
& +o(\alpha_0).
\end{aligned}
\end{equation*}

The goal is to determine $\nu_1, \nu_2$ such that the corresponding operators
$\Lcal_a$ and ${\tLa}^*$ satisfy \eqref{id1}.
Using the two expansions for $\Lcal_a$ and ${\tLa}^*$ it follows that
\begin{equation} \label{aux-1}
\begin{aligned}
{\tLa}^*\Lcal_a[\phi(s)](t) = & \phi(t) \\
&+ \alpha_0 \left\lbrace \Fourier^{-1}\left[ \left( \lambda_1(\omega) + \nu_1(-\omega)+\nu_2(-\omega)\right)
\Fourier[\phi(s)](\omega)\right](t) \right.  \\
& \qquad + \Fourier^{-1}\left[ -\imaginary \omega\lambda_1(\omega) \Fourier[s\phi(s)](\omega)\right] (t) \\
& \left. \qquad + t \Fourier^{-1}\left[\imaginary \omega \nu_1(-\omega) \Fourier[\phi(s)](\omega)\right](t)
\right\rbrace \\
&+o(\alpha_0)
\end{aligned} \end{equation}
To satisfy \eqref{id1} we require that the first order term in $\alpha_0$ of the equation \eqref{aux-1} has to vanish.
By taking the Fourier transform of this term we see that the term vanishes if 
\begin{equation} \label{aux-2}
\big( \lambda_1(\omega) + \nu_1(-\omega)+\nu_2(-\omega)\big)  \Fourier[\phi](\omega) - \omega\lambda_1(\omega) \frac{d
\Fourier[\phi]}{d\omega}(\omega)+
\frac{d}{d\omega} \left( \omega \nu_1(-\omega) \Fourier[\phi](\omega) \right)=0 ,
\end{equation}
where we have used the property \eqref{eq:diff} with n=1. Now, it is straightforward to see that the solution of the following system
\begin{equation} \label{cond-0}
\begin{aligned}
\lambda_1(\omega)+\nu_1(-\omega)+\nu_2(-\omega) + \frac{d (\omega \nu_1(-\omega))}{d\omega} &= 0\,,\\
-\omega\lambda_1(\omega)+\omega \nu_1(-\omega)&=0\,,
\end{aligned}
\end{equation}
satisfies the equation \eqref{aux-2} which directly implies \eqref{id1}. Equivalently, we get the following conditions 
\begin{equation} \label{cond-1}
\nu_1(-\omega) = \lambda_1(\omega) \text{ and } \nu_2(-\omega) = -3 \lambda_1(\omega) - \omega \frac{d
\lambda_1(\omega)}{d\omega}\;.
\end{equation}

Now, we introduce $\Gamma_\omega(x,y) $ and $\tilde\Gamma_\omega^a(x,y) $ which are the fundamental solutions of the
Helmholtz equations
\begin{equation} \label{fundamental-Helm} \omega^2 \Gamma_\omega(x,y)+\Delta_y \Gamma_\omega(x,y)= - \delta_x(y), \qquad
y\in \R^3 \end{equation}
and
\begin{equation} \label{fundamental-treq} \tilde{\kappa}(\omega)^2 \tilde\Gamma_\omega^a(x,y)+\Delta_y
\tilde\Gamma_\omega^a(x,y)= - \lambda(\omega) \delta_x(y), \qquad y\in \R^3, 
\end{equation}
respectively.
We can prove that 
\begin{equation}
\label{proposition-1} 
\pd{\tilde\Gamma^a}{t}=\tLa\left[ \pd{\Gamma}{t}\right],
\end{equation}
where $\tLa$ is defined in \eqref{tilde-La-operator}, 
\begin{equation} \label{Gamma}  \Gamma(x,y,t,\tau)=
\Fourier^{-1}\left\{ \Gamma_\omega(x,y) \right\}(t-\tau),\end{equation} and  \begin{equation} \label{Gamma^a}
\tilde\Gamma^a(x,y,t,\tau)= \Fourier^{-1}\left\{ \tilde\Gamma_\omega^a(x,y) \right\}(t-\tau),\end{equation} and
$\Fourier^{-1}$ denotes the inverse Fourier transform with respect to $\omega$. 

Then, we define the function $v_s^a(x,t)$ by
\begin{equation} \label{v_s^a}
v_s^a(x,t)=-\frac{1}{2\pi}\int_\R\int_{\partial\Omega} \imaginary \omega \tilde{\Gamma}_\omega^a(x,y) g_a(y,T-s)
d\sigma(y)  e^{-\imaginary \omega(t-s)} d\omega ,
\end{equation}
where we recall that $g_a(y,t):=p^a(y,t) \text{ for all } y\in\partial \Omega \text{ and } t\in[0,T).$

For the KSB model, it follows from \eqref{eq:kappa} and \eqref{order} that
\begin{equation*}
\lambda_1(\omega)=-\left(1+(-\imaginary \tau_0 \omega)^{\gamma-1}\right)^{-1/2}.
\end{equation*}
Using this in \eqref{cond-1} we get 
\begin{equation*}
\nu_1(\omega)=-\left(1+(\imaginary \tau_0 \omega)^{\gamma-1}\right)^{-1/2}.
\end{equation*}
and
\begin{equation*}
\nu_2(\omega) =\frac{7-\gamma}{2}\left(1+(\imaginary \tau_0 \omega)^{\gamma-1}\right)^{-1/2} +
\frac{\gamma-1}{2}\left(1+(\imaginary \tau_0 \omega)^{\gamma-1}\right)^{-3/2}\;.
\end{equation*}

Using these expressions for $\nu_1(\omega)$ and $\nu_2(\omega)$ we suggest the following choice for 
$\tilde{\kappa}(\omega)$ and $\lambda(\omega)$:
\begin{equation} 
\label{NSW-tilde-kappa} 
\tilde{\kappa}(\omega)=\omega (1 - \alpha_0 \nu_1(\omega)) \text{ and } 
\lambda(\omega)=1 + \alpha_0 \nu_2 (\omega)\,,
\end{equation} 
respectively. Note, that there is no remainder term in \eqref{cond-1}.

By applying the
expressions \eqref{NSW-tilde-kappa} into \eqref{fundamental-treq}, then the previous
definition of the function $v_s^a(x,t)$ yields the following identity:
\begin{equation*}
\begin{split}
& \left( \tilde{L}^{1/2} \left( \alpha_0 \Id + \tilde{L}^{1/2}\right)^2 \frac{\p^2}{\p t^2} -\tilde{L}^{1/2} \ \tilde{L}
\ \Delta \right) v_s^a(x,t)\\
 & = \bigg( \tilde{L}^{1/2} \ \tilde{L} +\alpha_0 \left( (\gamma-1) \Id +(7-\gamma) \tilde{L} \right) \bigg)
 \pd{\delta_s}{t}\left(g_a(x,T-s)\delta_{\partial \Omega}\right) \text{ for } x \in \Omega\;.
\end{split}
\end{equation*}
Here 
\begin{equation}
\label{eq:tildeL}
 \tilde{L}=\Id + (-\tau_0)^{\gamma-1}D_t^{\gamma-1}\,,
\end{equation}
and 
\begin{equation*}
 g_a(y,t):=p^a(y,t) \text{ for all } y\in\partial \Omega\,, t\in[0,T]\,,
\end{equation*}
where $T$ is supposed to be sufficiently large such that $p^a(x,t)=0=\pd{p^a}{t}(x,t)$ for $t\geq T$
and $x\in \Omega$. 

In the following subsection we prove that the functional $$ \Ical^a(x)=\int_0^T v_s^a(x,T)ds $$ is an approximation of
the initial datum $f(x)$. 
For doing this we have to make some regularization of the relevant operators.

\subsection{The regularized time reversal functional}
Matters of convergence of some infinite integrals defined above, suggest a regularization in the same way as in
\cite{AmmBreGarWah11} and \cite{Wah11}. We define the function 
\begin{equation} \label{v-sr^a}
 v_{s,\rho}^a(x,t)=-\frac{1}{2\pi}\int_{-\rho}^{\rho}\int_{\partial\Omega} \imaginary \omega
 \tilde{\Gamma}_\omega^a(x,y) g_a(y,T-s) d\sigma(y)
 e^{-\imaginary \omega(t-s)}
 d\omega .\end{equation}
as an approximation of $v_s^a(x,t)$ defined in \eqref{v_s^a}. Moreover, we define the regularized fundamental solution of $\tilde\Gamma^a$ as in \eqref{Gamma^a}
\begin{equation*}\tilde{\Gamma}_\rho^a(x,y,s,t)= \frac{1}{2\pi}\int_{-\rho}^{\rho} e^{-\imaginary \omega (t-s)}
\tilde{\Gamma}_\omega^a(x,y)d\omega ,\end{equation*}
the regularized operator $\widetilde\Lcal_{a,\rho}$ defined in \eqref{operator1}
\begin{equation*}
\widetilde\Lcal_{a,\rho}[\phi](t)=\frac{1}{2\pi} \int_0^\infty \phi(s)   \int_{-\rho}^\rho \frac{\omega
\lambda(\omega)}{\tilde\kappa(\omega)} e^{\imaginary\tilde\kappa(\omega)s}e^{-\imaginary\omega t}d\omega ds,
\end{equation*}
and its adjoint
\begin{equation*} \label{operator-adjoint}\widetilde\Lcal_{a,\rho}^*[\phi](t)=\frac{1}{2\pi} \int_{-\rho}^\rho
\frac{\omega \lambda(\omega)}{\tilde\kappa(\omega)} e^{\imaginary\tilde\kappa(\omega) t} \int_0^\infty
e^{-\imaginary\omega s}\phi(s)dsd\omega .\end{equation*}
Using these definitions we write the approximated version of the equations \eqref{id1} and \eqref{proposition-1},
respectively, that is

\begin{equation}\label{id-rho} 
\widetilde\Lcal_{a,\rho}^* \Lcal_a[\phi](t)= S_\rho[\phi](t) + o(\alpha_0)\,,
\end{equation} 
and
\begin{equation}
\label{proposition-rho} \pd{\tilde\Gamma^a_\rho}{t}=\widetilde\Lcal_{a,\rho}\left[ \pd{\Gamma}{t}\right] ,\end{equation}
where the operator $S_\rho\left[\phi\right]$ and the function $\Gamma$ were defined in \eqref{S_r} and \eqref{Gamma}, respectively. 

Similarly to the previous subsection, applying the definition \eqref{v-sr^a} of $v_{s,\rho}^a(x,t)$ in equation
\eqref{fundamental-treq} with usage of the expressions \eqref{NSW-tilde-kappa} and \eqref{NSW-tilde-kappa} we obtain the
following wave equation
\begin{equation}
\begin{split}
& \left( \tilde{L}^{1/2} \left( \alpha_0 \Id + \tilde{L}^{1/2}\right)^2 \frac{\p^2}{\p t^2} -\tilde{L}^{1/2} \ \tilde{L}
\ \Delta \right) v_{s,\rho}^a(x,t)\\
 & = S_\rho\left[ \bigg( \tilde{L}^{1/2} \ \tilde{L} +\alpha_0 \left( (\gamma-1) \Id +(7-\gamma) \tilde{L} \right)
 \bigg) \pd{\delta_s}{t}\right] \left(g_a(x,T-s)\delta_{\partial \Omega}\right) \text{ for }x\in \Omega\;.
\end{split}
\end{equation}
Finally, we can obtain the reconstruction functional $\Ical_\rho^a$. Indeed, since equations \eqref{id-rho} and
\eqref{proposition-rho} hold, then Proposition 2.3.5 in \cite{Wah11} shows that
\begin{equation*}
\Ical_\rho^a(x)=\int_0^T v_{s,\rho}^a(x,T)ds \longrightarrow f(x), \qquad \text{as }
\rho\to+\infty.
\end{equation*}

\begin{remark} The latter proposition suggest that the larger $\rho$ we choose, the better approximation we get.
However, the previous computation of $\Ical_\rho^a(x)$ involves the integration of the fundamental solution 
$\tilde{\Gamma}_\omega^a(x,y)$ which grows exponentially as $\exp\left\lbrace \Im\{\tilde{\kappa}(\omega)\}|x-y|\right\rbrace  $. In order to ensure stability of $\Ical_\rho^a(x)$ this term must not be greater than one \cite[Remark 2.3.6]{Wah11}. For large $\omega$, the expression $| \Im\{\tilde{\kappa}(\omega)\} |$ is behaving like (and is less than)  $\alpha_0
|\omega|\sin\frac{(\gamma-1)\pi}{4}, \ \gamma \in (1,2].$ So one should not use frequencies larger than
$\frac{1}{\alpha_0 diam(\Omega)}$. Hence, we get $\rho \simeq\frac{1}{\alpha_0 diam(\Omega)}$ to be the threshold for the imaging
functional stability, where $diam(\Omega)$ denotes the diameter of the domain $\Omega$. A finer estimation of the
threshold can be given if we use that $| \Im\{\tilde{\kappa}(\omega)\} |$ is behaving like $\alpha_0
\tau_0^{\frac{1-\ga}{2}} |\omega|^{\frac{3-\ga}{2}}
\sin\frac{(\gamma-1)\pi}{4}, \ \gamma \in (1,2]$, for large values of $\omega$. Consequently, \begin{equation*}\rho
\simeq\dfrac{\tau_0^{\frac{\ga-1}{3-\ga}}}{\left( \alpha_0 diam(\Omega) \sin\frac{(\gamma-1)\pi}{4}\right)^{\frac{2}{3-\ga}}}
,\end{equation*} with $\gamma \in (1,2]$.\end{remark}

\section{The NSW model} 
Let $p^a$ satisfy the following problem
\begin{equation} \label{nsw-wave_eq}
\left(\Id + \tilde{\tau} \frac{\p }{\p t} \right) \frac{\p^2 p^a}{\p t^2} -\left(\Id + \tau \frac{\p}{\p t} \right)
\Delta p^a = \left(\Id + \tau \frac{\p}{\p t} \right)\pd{\delta_0}{t} f,
\end{equation}
along with the conditions in \eqref{wave-KSB}, where we consider the NSW model for one relaxation process, as defined in
\cite{KowSch10}. Moreover, we assume that $\tau>\tilde\tau>0$ so that the strong causality condition in \cite{KowSch10} is satisfied
and that $\tau$ and $\tilde\tau$ are small and of the same magnitude.

Substituting $\tilde\tau = \alpha_0 \tilde r$ and $\tau = \alpha_0 r$, we find 
\begin{equation*}
\kappa(\omega)=\omega \sqrt{\dfrac{1-\imaginary \omega \tilde\tau}{1-\imaginary \omega \tau}} = \omega \sqrt{\dfrac{1-\imaginary \alpha_0 \omega \tilde r}{1-\imaginary \alpha_0 \omega r}}.
\end{equation*}
By applying the expansion \eqref{order} for $\kappa(\omega)$ in terms of $\alpha_0$ we get $$\lambda_1(\omega)=-\imaginary \omega \frac{r-\tilde r}{2}.$$
Consequently, the auxiliary functions $\nu_i,\  i=1,2$  are given by the conditions \eqref{cond-1} and read as follows
\begin{equation} \label{nsw-nu_1,2}
\nu_1(\omega)=\imaginary \omega \frac{r-\tilde r}{2} \ \text{ and } \
\nu_2(\omega)=-2\imaginary\omega (r-\tilde r)\;.
\end{equation}
Therefore, we obtain the expansion of $\tilde{\kappa}(\omega)$ and $\lambda(\omega)$ from \eqref{tkappa}, i.e.,
\begin{equation*}
\begin{aligned}
\tilde{\kappa}(\omega) &:=\omega(1-\alpha_0 \nu_1(\omega))+{\mathcal O}(\alpha_0^2)\,,\\
\lambda(\omega)&:=1+\alpha_0 \nu_2(\omega)+{\mathcal O}(\alpha_0^2),
\end{aligned}
\end{equation*}
 where now the functions $\nu_i,\  i=1,2$ are given by \eqref{nsw-nu_1,2}.

We make the following choices for the functions $\tilde{\kappa}(\omega)$ and $\lambda(\omega)$
\begin{equation} \label{nsw-tka+la}
\tilde{\kappa}(\omega)=\omega \sqrt{\dfrac{1+\imaginary \omega \tilde\tau}{1+\imaginary \omega \tau}} \ \text{ and } \
\lambda(\omega)=\left( \dfrac{1+\imaginary \omega \tilde\tau}{1+\imaginary \omega \tau}\right) ^2,
\end{equation}
which satisfy \eqref{tkappa} for non-vanishing ${\mathcal O}(\alpha_0^2)$ terms.
From here, using the Helmholtz equation \eqref{fundamental-treq} and with the same arguments as in the previous section
we derive the following regularized time reverted attenuated equation (which is of course not unique)
\begin{equation}
\label{nsw-reg-treq}
\begin{aligned}
~ & \left( \left(\Id - \tau \frac{\p }{\p t} \right)\left(\Id - \tilde{\tau} \frac{\p }{\p t} \right) \frac{\p^2 }{\p t^2}
-\left(\Id - \tau \frac{\p}{\p t} \right)^2 \Delta \right) v_{s,\rho}^a(x,t) \\
= & S_\rho\left[ \left(\Id - \tilde\tau
\frac{\p}{\p t} \right)^2 \pd{\delta_s}{t}\right] \left(g_a(x,T-s)\delta_{\partial \Omega}\right),
\end{aligned}
\end{equation}
where $S_\rho$ is defined in \eqref{S_r}. Here \begin{equation*}
 g_a(y,t):=p^a(y,t) \text{ for all } y\in\partial \Omega\,, t\in[0,T]\,,
\end{equation*}
where $p^a(y,t)$ satisfies equation \eqref{nsw-wave_eq} and $T$ is supposed to be sufficiently large such that $p^a(x,t)=0=\pd{p^a}{t}(x,t)$ for $t\geq T$
and $x\in \Omega$. 

The reconstruction imaging functional is given by \eqref{functional-reg}, i.e.
\begin{equation*}
\Ical_\rho^a(x)=\int_0^T v_{s,\rho}^a(x,T)ds \longrightarrow f(x), \qquad \text{as } \rho\to+\infty\,.
\end{equation*}

\begin{remark}
For the case of N relaxation processes, the procedure will be conceptually the same. Now, the attenuated wave equation
has a more complicated form \cite{NacSmiWaa90} and we have 
\begin{equation*}
\kappa(\omega) = \omega \sqrt{\dfrac{1}{N}\sum_{j=1}^N\dfrac{1-\imaginary \omega \tilde\tau_j}{1-\imaginary \omega \tau_j}}.
\end{equation*}
Here, we assume that $\{\tau_j,\tilde{\tau}_j\}_1^N$ are small 
and of the same magnitude, i.e. all  $\{\tau_j,\tilde{\tau}_j\}_1^N$ are of order ${\mathcal O}(\alpha_0)$. 
Then the expansions of the functions $\kappa(\omega)$, $\tilde{\kappa}(\omega)$ and $\lambda(\omega)$ in terms of $\alpha_0$, as they were given in \eqref{order} and \eqref{tkappa}, 
and the relevant asymptotic analysis allow us to make (as previously) the following choices
\begin{equation} \label{nsw-tka+la-gen}
\tilde{\kappa}(\omega)=\omega \sqrt{\dfrac{1}{N}\sum_{j=1}^N\dfrac{1+\imaginary \omega \tilde\tau_j}{1+\imaginary \omega \tau_j}} 
\text{ and } \lambda(\omega)=\left( \dfrac{1}{N}\sum_{j=1}^N\dfrac{1+\imaginary \omega \tilde\tau_j}{1+\imaginary \omega
\tau_j}\right)^2, \end{equation} 
which lead to the corresponding time-reverted attenuated equation. One can find the corresponding wave equation by applying the 
inverse Fourier transform on the Helmholtz equation \eqref{fundamental-treq}, for the later values of 
$\tilde{\kappa}(\omega)$ and $\lambda(\omega)$. This procedure will lead to an equation similar to \eqref{nsw-reg-treq},
but this time one obtains a more complicated form.

Following the analysis of the previous section we will choose the value of the truncation parameter $\rho$ appearing in \eqref{nsw-reg-treq} by finding
the behaviour of the expression $|\Im (\tilde{\kappa}(\omega))|$, with $\tilde{\kappa}(\omega)$ given in \eqref{nsw-tka+la}. For $|\omega| \tilde\tau < |\omega|\tau<1$ we find that
$|\Im(\tilde{\kappa}(\omega))|$ is behaving like $\dfrac{\omega^2}{2}(\tau-\tilde\tau)$, which (following the arguments
from \cite[Remark 2.3.6]{Wah11}) yields the truncation parameter
\begin{equation*}
\rho\simeq\frac{1}{\sqrt{diam(\Omega)(\tau-\tilde\tau)}}.
\end{equation*}
 Moreover, this value of the truncation parameter is sufficient
for the other two cases, i.e. for $1<|\omega|\tilde\tau<|\omega|\tau$ and $|\omega|\tilde\tau<1<|\omega|\tau$, when
$|\Im(\tilde{\kappa}(\omega))|$ is behaving like 
$\frac{1}{2}\sqrt{\frac{\tilde\tau}{\tau}}\left(\frac{1}{\tilde\tau}
-\frac{1}{\tau}\right) $ and $\sqrt{\dfrac{|\omega|}{2\tau}}\left| -1+\frac{1}{2}\left( |\omega|\tilde\tau +
\frac{1}{|\omega|\tau}\right)\right|$, respectively.

For the case of N relaxation processes, one can use similar arguments to get an estimation of the truncation parameter
$\rho$. In the case that $|\omega|\tilde\tau_j<|\omega|\tau_j<1,$ for all $j=1,\ldots,N$ we find that
$|\Im(\tilde{\kappa}(\omega)\}|$, with $\tilde{\kappa}(\omega)$ given in \eqref{nsw-tka+la-gen}, is behaving like $\dfrac{\omega^2}{2N}\sum_{j=1}^N(\tau_j-\tilde\tau_j)$, which yields
the truncation parameter
\begin{equation*}
\ds\rho\simeq\sqrt{\dfrac{N}{diam(\Omega)\sum_{j=1}^N(\tau_j-\tilde\tau_j)}}.
\end{equation*} 
For the several other cases, the arguments of the previous remark, along with the usage of the
triangular inequality, give us the opportunity to observe that the above-mentioned estimation for the truncation
parameter is sufficient.
\end{remark}

\begin{remark}
The formal procedure outlined above also applies to thermo-viscous model, as this was defined by the wave equation (84) in \cite{KowSch10}. 
This is the special case of the NSW-model with one relaxation process and $\tilde\tau=0$. Note that the thermo-viscous wave equation (84) in \cite{KowSch10} refers to 
a not strongly causal model and has a different RHS from the thermo-viscous equation \eqref{wave-eq-att-tv}, which was analysed in \cite{Wah11}. 
\end{remark}

\section{Higher order terms}

In this section we describe the procedure for evaluating higher order terms of the operators $\Lcal_a$ and $\tLa^*$, defined in \eqref{operator1} and \eqref{tilde-La-operator} and 
consequently a higher order approximation of the reconstruction functional $\Ical_\rho^a$.

In particular, this method allows determining the higher order terms of the asymptotic expansion for the functions 
$\tilde{\kappa}(\omega)$ and $\lambda(\omega)$, appearing in \eqref{tilde-La-operator}.
We make the following ansatz
\begin{equation}
\label{ka-gen}\kappa(\omega)=\omega \sum_{j=0}^\infty (-1)^j \lambda_j(\omega) a^j, \qquad \lambda_0(\omega)=1
\end{equation}
and consequently
\begin{equation}
\label{w/ka-gen}\frac{\omega}{\kappa(\omega)}=\sum_{j=0}^\infty  \mu_j(\omega) a^j,
\end{equation}
with
\begin{equation}
\label{mu0,1,2-gen} \mu_0(\omega)=1, \qquad \mu_1(\omega)=\lambda_1(\omega), \qquad
\mu_2(\omega)=\lambda_1^2(\omega)-\lambda_2(\omega)
\end{equation}
and in general $$\mu_j(\omega)=\Pcal_j\left( \left\{ \lambda_k(\omega) \right\}_{k=1}^j \right), $$ where $\Pcal_j$ denotes
a polynomial (on several variables) of order $j$. In addition, the expansion \eqref{ka-gen} yields the following
\begin{equation}
\label{exp-ka-gen}e^{\imaginary \kappa(\omega) s}=e^{\imaginary\omega s}\sum_{j=0}^\infty \psi_j(\omega,s) a^j,
\end{equation}
with
\begin{equation}
\label{psi0,1,2-gen} \psi_0(\omega,s)=1, \qquad \psi_1(\omega,s)=-\imaginary\omega\lambda_1(\omega)s, \qquad
\psi_2(\omega,s)=\imaginary\omega\lambda_2(\omega)s+\frac{1}{2}\left( \imaginary\omega\lambda_1(\omega)s\right)^2 
\end{equation}
and in general 
\begin{equation}
\label{psij-gen}\psi_j(\omega,s)=\Qcal_j\left( \left\{ \imaginary\omega\lambda_k(\omega) s \right\}_{k=1}^j \right),
\end{equation} where $\Qcal_j$ denotes a polynomial (on several variables) of order $j$.
In the case of higher order terms we make the ansatz
\begin{equation}
\label{wka-gen}\tilde\kappa(\omega)=\omega \sum_{j=0}^\infty (-1)^j \lambda_j(-\omega) a^j.
\end{equation}
This expansion is without loss of generality. In the procedure described in section \ref{section-KSB}, concerning the first order approximation of $\tilde\kappa(\omega)$, 
the relevant terms were considered unknown and one had to determine them. However, applying the previous expansion we
state a consistency condition for the time-reversal algorithm which, instead of giving a system of equations, yields a
set of identities, as we will see below.

In the same way with the previously mentioned expansions we get 
\begin{equation}
\label{w/wka-gen}\frac{\omega}{\tilde\kappa(\omega)}=\sum_{j=0}^\infty \mu_j(-\omega) a^j
\end{equation}
and
\begin{equation}
\label{exp-wka-gen}e^{\imaginary \tilde\kappa(\omega) s}=e^{\imaginary\omega s}\sum_{j=0}^\infty
\widetilde\psi_j(\omega,s) a^j,
\end{equation}
with
\begin{equation}
\label{wpsij-gen} \widetilde\psi_j(\omega,s)=\Qcal_j\left( \left\{ \imaginary\omega\lambda_k(-\omega) s \right\}_{k=1}^j \right), 
\end{equation}
where $\Qcal_j$ denote the same polynomials as in \eqref{psij-gen} .

Finally, we consider the expansions
\begin{equation}
\label{la-gen}\lambda(\omega)=\sum_{j=0}^\infty \beta_j(\omega) a^j, \qquad \beta_0(\omega)=1
\end{equation}
and
\begin{equation}
\label{gamma-gen}\frac{\omega}{\tilde\kappa(\omega)}\lambda(\omega)=\omega \sum_{j=0}^\infty \ga_j(\omega) a^j, \qquad
\ga_0(\omega)=1.
\end{equation}
These, expansions along with \eqref{wka-gen} yield the following relation
\begin{equation}
\label{beta-gen}\beta_n(\omega)=\sum_{i+j=n} (-1)^j \ga_i(\omega) \lambda_j(-\omega).
\end{equation}

Now, since we know explicitly $\tilde\kappa(\omega)$, our strategy consists of determining $\lambda(\omega)$; consequently
the Helmholtz equation \eqref{fundamental-treq}, under the procedure described in the relevant section, provides the corresponding time-reverted wave equation.
So, in what follows out target is to describe a procedure to find $\beta_n(\omega), \ n\in \N$, equivalently determine
the terms of the asymptotic expansion of $\lambda(\omega)$.

Applying the previous expansions in the expressions \eqref{operator1} and \eqref{operator2}, we get 
\begin{equation}
\label{op-gen} \Lcal_a[\phi](t)=\sum_{k=0}^\infty f_k[\phi](t) a^k  \text{ and }
\tLa^*[\phi](t)=\sum_{k=0}^\infty g_k[\phi](t) a^k,
\end{equation} 
with $f_0 \equiv g_0 \equiv Id$, $\{f_k,g_k\}_{k=1}^n$ being operators that can be obtained explicitly in terms of
$\{\lambda_j,\mu_j,\psi_j,\widetilde\psi_j,\ga_j\}_{j=1}^k$, for $k=1,\ldots,n,$ respectively. In particular,
$\{f_k,g_k\}_{k=1}^2$ are obtained explicitly with use of the expressions \eqref{mu0,1,2-gen}, \eqref{psi0,1,2-gen},
\eqref{wpsij-gen}, i.e., 
$$f_1[\phi(s)](t)=\Fourier^{-1}\left[ -\imaginary\omega\lambda_1(\omega)
\Fourier[s\phi(s)](\omega)+\lambda_1(\omega)\Fourier[\phi(s)](\omega)\right](t),$$
\begin{equation} \begin{split}f_2[\phi(s)](t)=\Fourier^{-1}\bigg[\left( \lambda_1^2(\omega)-\lambda_2(\omega) \right)
\Fourier[\phi(s)](\omega) -\imaginary \omega \left( \lambda_1^2(\omega)-\lambda_2(\omega) \right) \Fourier[s\phi(s)](\omega) \\
+\frac{1}{2}\left( \imaginary\omega\lambda_1(\omega)\right)
^2\Fourier[s^2\phi(s)](\omega)\bigg](t),\end{split}\nonumber\end{equation}
$$g_1[\phi(s)](t)=\Fourier^{-1}\left[ \ga_1(-\omega) \Fourier[\phi(s)](\omega)\right] + t \Fourier^{-1}\left[
\imaginary\omega\lambda_1(\omega)\Fourier[\phi(s)](\omega)\right](t),$$
\begin{equation} \begin{split}g_2[\phi(s)](t)=\Fourier^{-1}\left[ \ga_2(-\omega) \Fourier[\phi(s)](\omega)\right] + t
\Fourier^{-1}\left[ \imaginary\omega \left( \ga_1(-\omega)\lambda_1(\omega)- \lambda_2(\omega)\right)
\Fourier[\phi(s)](\omega)\right](t) \\ + t^2 \Fourier^{-1}\left[ \frac{1}{2}\left( \imaginary\omega\lambda_1(\omega)\right)^2
\Fourier[\phi(s)](\omega)\right](t).\end{split}\nonumber\end{equation}
Consequently, we get the following asymptotic expansion
\begin{equation}
\label{op-gen-id} \tLa^*\Lcal_a[\phi](t)=\phi[t] +\bigg(f_1[\phi](t)+ g_1[\phi](t)\bigg) a+
\bigg(f_2[\phi](t)+ g_2[\phi](t)+g_1[f_1[\phi]](t)\bigg) a^2 + o(a^2),
\end{equation} 
when $a\to 0$, with $\{f_k,g_k\}_{k=1}^2$ defined above. The identity
\begin{equation}
\label{id-gen} \tLa^*\Lcal_a[\phi](t)=\phi(t) + o(a^n), \qquad a\to 0
\end{equation} 
is satisfied, up to the 1st order (equivalently the identity \eqref{id1} is satisfied), when
$$\Fourier\left[ f_1[\phi] +  g_1[\phi]\right](\omega)=0, \quad \forall \omega \in \R.$$ Using the definitions of the
operators $f_1$ and $g_1$ given above and the property \eqref{eq:diff} this condition yields the following equation
\begin{equation} \label{cond-gen-1} \ga_1(-\omega)=-2\lambda_1(\omega)-\omega\frac{d\lambda_1(\omega)}{d\omega},
\end{equation}
which is equivalent with the conditions \eqref{cond-0}. 
By applying the condition \eqref{cond-gen-1} in the expression \eqref{beta-gen} for $n=1$, we get that
\begin{equation}\label{beta_1}\beta_1(-\omega)=-\lambda_1(\omega)+\gamma_1(-\omega)=-3\lambda_1(\omega)-\omega\frac{d\lambda_1(\omega)}{d\omega},\end{equation}
which is the first order approximation, of the function $\lambda(\omega)$, as this was defined in \eqref{la-gen}. From the expressions \eqref{tkappa} and \eqref{la-gen}, one can see that $\nu_2\equiv \beta_1$, as they were introduced in
these equalities, respectively. So, equation \eqref{beta_1} is exactly the second of the conditions \eqref{cond-1}. Note
that the first of these conditions, is satisfied as an identity by the choice of the expansion of
$\tilde\kappa(\omega)$, given in
\eqref{wka-gen}.

Following the same procedure for the 2nd order in the identity \eqref{id-gen}, the expression \eqref{op-gen-id} yields
$$\Fourier [ f_2[\phi]+ g_2[\phi]+g_1[f_1[\phi]]](\omega)=0,$$ which after some calculations where we made use of the
definitions of the operators $f_2$, $g_2$ and the fact that
\begin{equation} \begin{split} g_1[f_1[\phi(s)]](t)=\Fourier^{-1}\big[ \ga_1(-\omega)\lambda_1(\omega)
\Fourier[\phi(s)](\omega)-\imaginary\omega\ga_1(-\omega)\lambda_1(\omega) \Fourier[s\phi(s)](\omega)\big] \\ + t
\Fourier^{-1}\left[ \imaginary\omega  \lambda_1^2(\omega)\Fourier[\phi(s)](\omega)- \left( \imaginary\omega\lambda_1(\omega)\right)^2
\Fourier[s\phi(s)](\omega)\right](t) 
\end{split} \nonumber\end{equation}
we obtain the following equation 
\begin{equation} \label{cond-gen-2}
\begin{split}\ga_2(-\omega)&=-\lambda_1^2(\omega)+\lambda_2(\omega)-\ga_1(-\omega)\lambda_1(\omega)\\ 
&+\frac{d}{d\omega}\bigg(\omega \left( -\lambda_1^2(\omega)+\lambda_2(\omega)-\ga_1(-\omega)\lambda_1(\omega)\right)
\bigg) -\frac{1}{2}\frac{d^2}{d\omega^2}\left( \big(\omega\lambda_1(\omega) \big) ^2\right).
\end{split}\end{equation}
Substitution of the conditions \eqref{cond-gen-1} and
\eqref{cond-gen-2} in the expression \eqref{beta-gen} provides the second order approximation for the function
$\lambda(\omega)$, i.e. $\beta_2(\omega)$, as this was defined in \eqref{la-gen}.

A general treatment of this problem would appear as follows: The general form of \eqref{op-gen-id} would be
$$ \tLa^* \Lcal_a[\phi](t) = \phi(t) +\sum_{k=1}^n h_k[\phi](t) a^k + o(a^n) ,  \qquad a \to 0 , $$ 
where $\{h_k\}_{k=1}^n$ are operators expressed explicitly (after some calculations and the use of the previous expansions) in terms of $\{\lambda_j,\ga_j\}_{j=1}^k$, for $k=1,\ldots,n$, respectively. Consequently, the condition
$$\Fourier[h_k[\phi]](\omega)=0, \qquad \forall \omega \in \R $$ yields $\ga_k(-\omega)$ as an
explicit expression of $\left\{ \lambda_j(\omega)\right\}_{j=1}^k$, $\left\{ \ga_j(-\omega)\right\}_{j=1}^{k-1}$ and
their derivatives up to order $k$. Hence, the expression \eqref{beta-gen} yields the  $k$-th order approximation for the
function $\lambda(\omega)$, i.e. $\beta_k(\omega)$.

\def\cprime{$'$} \providecommand{\noopsort}[1]{}

%
%

\section*{Acknowledgements}
\label{sec:acknowledgements}

This work has been supported by the Austrian Science Fund (FWF) within the national research networks Photoacoustic
Imaging in Biology and Medicine (project S10505) and Geometry and Simulation (project S11704).
\section*{Appendix}
In this paper we use the basic notation:

\begin{itemize}
 \item $\de_s$ denotes the one-dimensional Dirac distribution at $s \in \R$. Moreover, we denote by $\delta_x$ the 3-dimensional $\delta$-distribution with center $x$ - 
       $x$ is always a 3D variable.
 \item $\Fourier$ denotes the Fourier transform 
       \begin{equation}\label{def:Fourier}
       \Fourier[\phi](\omega)= \frac{1}{\sqrt{2\pi}}\int_\R e^{\imaginary \omega t} \phi(t)dt\;.
       \end{equation}
 \item If we want to specify the argument, with which respect the Fourier transform is applied it is specified as follows:
       \[\ds \Fourier[\phi(s)](\omega)=\frac{1}{\sqrt{2\pi}}\int_\R e^{\imaginary \omega s} \phi(s)ds.\]
       In this case it is with respect to the variable $s$.
\item \begin{equation}\label{S_r} 
      S_\rho\left[\phi\right](t)= \frac{1}{\sqrt{2\pi}}\int_{-\rho}^{\rho} e^{-\imaginary \omega t} \Fourier[\phi](\omega)d\omega\;.
      \end{equation}
\item $\delta_{\partial \Omega}$ denotes the surface Dirac mass on $\partial \Omega$.
\end{itemize}
Moreover we use the following property of the Fourier transform
\begin{equation}
\label{eq:diff}
\Fourier[t^n \phi(t)](\omega)=(-\imaginary)^n \frac{d^n}{d\omega^n}\Fourier[\phi(t)](\omega),
\end{equation}
for $n \in \N$.

\end{document}